\documentstyle[12pt]{article}
\title{The Matching Number and Hamiltonicity of Graphs}
\author {Rao Li \\                
         Dept. of mathematical sciences \\
         University of South Carolina Aiken \\
	     Aiken, SC 29801 \\
         {\it Email: raol@usca.edu }
         }
\date{Feb. 25, 2017}

\begin{document}
\maketitle
\begin{abstract}
The matching number of a graph G is the size of a maximum matching in the graph.
 In this note, we present a sufficient condition involving the matching number for the Hamiltonicity of graphs.   
 \end{abstract} 
$$2010 \,\, Mathematics \,\, Subject \,\, Classification: \,\, 05C70, \, 05C45.$$
$$Keywords: Matching \,\, Number, \,\, Hamiltonicity$$ \\

We consider only finite undirected graphs without loops or multiple edges.
Notation and terminology not defined here follow those in \cite{Bondy}.
Let $G = (V(G), \, E(G))$ be a graph. A matching $M$ in $G$ is a set of pairwise non-adjacent edges. 
A maximum matching is a matching that contains the largest possible number of edges. 
The matching number, denoted $m(G)$, of a graph $G$ is the size of a maximum matching.
For a vertex $u$ and a vertex subset $U$ in $G$, we use $N_U(u)$ to denote all the neighbors of $u$ in $U$.
We use $G_1 \vee G_2$ to denote the the join of two disjoint graphs $G_1$ and $G_2$. 
We define $\cal{F}$ $ := \{ G : K_{p, \, q} \subseteq G \subseteq K_p \vee (qK_1),$ where $q \geq p + 1 \geq 3\}$. 
 A cycle $C$  in a graph $G$ is called a Hamiltonian cycle of $G$ if $C$ contains all the vertices of $G$.  
A graph $G$ is called Hamiltonian if $G$ has a Hamiltonian cycle. \\

The purpose of this note is to present a sufficient condition based on 
the matching number for the Hamiltonicity of graphs. The main result is as follows.\\

\noindent {\bf Theorem 1.} Let $G$ be a graph of order $n \geq 3$ with matching number $m$ 
and connectivity $\kappa$ ($\kappa \geq 2$). If
$m \leq \kappa$, then $G$ is Hamiltonian or $G \in \cal{F}$.\\

\noindent{\bf Proof of Theorem 1.} Let $G$ be a graph satisfying the conditions in Theorem $1$. 
Suppose $G$ is not Hamiltonian.
 Since $\kappa \geq 2$, $G$ contains a cycle. 
Choose a longest cycle $C$ in $G$ and give an orientation on $C$. For a vertex $u$ on $C$, we use
$u^+$ to denote the successor of $u$ along the direction of $C$. $u^{+2}$ is defined as 
the successor of $u^+$ along the direction of $C$. 
Since $G$ is not Hamiltonian, 
there exists a vertex $x_0 \in V(G) \backslash V(C)$. By Menger's theorem, 
we can find $s$ ($s \geq \kappa$) pairwise disjoint (except for $x_0$) paths $P_1$, $P_2$, ..., $P_s$ 
between $x_0$ and $V(C)$. Let $u_i$ be the end vertex of $P_i$ on $C$, where $1 \leq i \leq s$.  
Then a standard proof in Hamiltonian graph theory yields that $S := \{x_0, u_1^+, u_2^+, ..., u_s^+ \}$ 
is independent (otherwise $G$ would have cycles which are longer than $C$).
Obviously, the edges of $u_1u_1^+$, $u_2u_2^+$, ..., and $u_su_s^+$ form a matching in $G$. Thus
$\kappa \leq s \leq m \leq \kappa$. Therefore $\kappa = s = m$. The remainder of proofs consists of 
five claims and their proofs. \\

\noindent {\bf Claim 1.} Let $H$ be the component in $V(G) \backslash V(C)$ that contains $x_0$. 
Then $H$ consists of the singleton $x_0$. \\

\noindent {\bf Proof of Claim 1.} Suppose, to the contrary, that Claim $1$ is not true. Then we can find an edge, say $e$, in $H$. Then 
the edges of $e$, $u_1u_1^+$, $u_2u_2^+$, ..., and $u_su_s^+$ form a matching in $G$, giving a contradiction of
$\kappa + 1 = s + 1 \leq m = \kappa$.  \\

\noindent {\bf Claim 2.}  $u_{i + 1} = u_i^{+2}$ for each $i$ with $1 \leq i \leq s$, where $u_{s + 1}$ is regarded as $u_1$.  
Namely, $C = u_1u_1^+u_2u_2^+ \, ... \, u_su_s^+u_1$.  \\

\noindent {\bf Proof of Claim 2.} Suppose, to the contrary, that there exists one $i$ with $1 \leq i \leq s$ such that 
$u_{i + 1} \neq u_i^{+2}$. Without loss of generality, we assume that $u_{2} \neq u_1^{+2}$. 
 Then
the edges of $x_0u_1$, $u_1^+u_1^{+2}$, $u_2u_2^+$, $u_3u_3^+$, ..., and $u_su_s^+$ 
form a matching in $G$, giving a contradiction of
$\kappa + 1 = s + 1 \leq m = \kappa$. \\

\noindent {\bf Claim 3.} If $V(G) \backslash (V(C) \cup \{\, x_0\, \})$ is not empty, then
$V(G) \backslash (V(C) \cup \{\, x_0\, \})$ is an independent set. \\

\noindent {\bf Proof of Claim 3.} Using the similar arguments as the ones in the proofs of Claim $1$, we can prove that
Claim $3$ is true. \\ 

\noindent {\bf Claim 4.} If the independent set $V(G) \backslash (V(C) \cup \{\, x_0\, \} 
:= \{\, x_1, x_2, ..., x_r\,\}$ is nonempty,
then $N_C(x_i) = \{\, u_1, u_2, ..., u_s\,\}$ for each $i$ with $1 \leq i \leq r$. \\

\noindent {\bf Proof of Claim 4.} Suppose, to the contrary, that there exists one $i$ with $1 \leq i \leq s$ such that 
$N_C(x_i) \neq \{\, u_1, u_2, ..., u_s\,\}$. Without loss of generality, we assume that $N_C(x_1) \neq \{\, u_1, u_2, ..., u_s\,\}$.
Using the similar arguments as the ones in the proofs of Claim $2$, we can prove that
$C = z_1z_1^+z_2z_2^+ \, ... \, z_sz_s^+z_1$, where $N_C(x_1) = \{\, z_1, z_2, ..., z_s\,\}$.
Since $N_C(x_1) \neq \{\, u_1, u_2, ..., u_s\,\}$, we must have that 
$N_C(x_1) = \{\, z_1, z_2, ..., z_s\,\} = \{\, u_1^+, u_2^+, ..., u_s^+\,\}$, In this case, we can easily find a cycle in
$G$ which is longer than $C$, giving a contradiction. \\

\noindent {\bf Claim 5.}  $u_i^+u_j \in E$ for each $i$ with $1 \leq i \leq s$, where $u_{s + 1}$ is regarded as $u_1$.    \\

\noindent {\bf Proof of Claim 5.}  If $i = 1$, it is obvious that $u_1^+u_1 \in E$ and $u_1^+u_2 \in E$. Suppose, to the contrary, that
there exists one $j$ with $3 \leq j \leq s$ such that $u_1^+u_j \not \in E$. 
Then $G[V(G) \backslash \{\, u_1, u_2, ..., u_{j - 1}, u_{j + 1}, ..., u_s\,\}]$ is disconnected, contradiction to the assumption that
the connectivity of $G$ is $\kappa$. Simiarly, we can prove that $u_i^+u_j \in E$ for each $i$ with $2 \leq i \leq s$, 
where $u_{s + 1}$ is regarded as $u_1$. \\

Claims $1$ - $5$ imply that $G \in \cal{F}$. So we complete the proof of Theorem $1$. \hfill{$\Box$} \\


\begin{thebibliography}{20}
  
\bibitem{Bondy} J. A. Bondy and U. S. R. Murty, Graph Theory with Applications, 
                Macmillan, London and Elsevier, New York (1976).
       
\end{thebibliography}
\end{document}